\documentclass{amsart}

\usepackage{amssymb}
\usepackage{amsmath}
\usepackage{amsfonts}
\usepackage{amsxtra}
\usepackage{amscd}
\usepackage{amsthm}
\usepackage{epsfig}
\usepackage{enumerate}
\usepackage{float}
\usepackage{xy}

\xyoption{all}

\theoremstyle{plain} \newtheorem{Thm}{Theorem}[section]
\theoremstyle{plain} 
\theoremstyle{plain} 
\theoremstyle{plain} 
\theoremstyle{definition} 
\theoremstyle{definition} \newtheorem{Exa}[Thm]{Example}
\theoremstyle{remark}

\newcommand{\End}{\mbox{$\mathrm{End}$}}

\newcommand{\newfontobj}[2]{
  \newcommand{#1}[1]{
    \expandafter\def\csname##1\endcsname{{#2 ##1}}}}

\newfontobj{\class}{\rm}
\newfontobj{\lang}{\bf}
\newfontobj{\oper}{\rm}

\class{PSPACE}		% Language classes have to implement
\class{AM}		% a clean way of doing the complement classes (simple)
\class{MA}
\class{NP}
\class{UP}
\class{P}
\class{RP}
\class{BPP}
\class{DTIME}
\class{ZPP}
\class{EXPSPACE}
\class{coNP}
\class{coRP}
\class{coAM}
\class{PH}
\class{co}
\class{BP}	

\lang{HN} 		
\lang{SAT}
\lang{PIM}
\lang{CSG}
\lang{ESC}
\lang{IP}
\lang{PRIMES}

\oper{Pr}    % These are the operators
\oper{E}     % Expectation
\oper{Ord}   % Order of an element/group
\oper{li}

\floatstyle{ruled}
\newfloat{Algorithm}{thp}{loa}[section]

\usepackage{fullpage}

\title{On the existence of distortion maps on ordinary elliptic curves}
\author{Denis Charles}

\begin{document}
\maketitle
\section{Introduction}
\parindent=0pt

An important problem in cryptography is the so called Decision Diffie-Hellman problem (henceforth abbreviated DDH). 
The problem is to distinguish
triples of the form $(g^a,g^b,g^{ab})$ from arbitrary triples from a cyclic group $G = \langle g \rangle$. It turns out that for (cyclic subgroups of) the group of $m$-torsion points
on an elliptic curve over a finite field, the DDH problem admits an efficient solution if there exists a suitable endomorphism called a distortion map 
(which can be efficiently computed) on
the elliptic curve.\\

Suppose $m$ is relatively prime to the characteristic of a finite field $\mathbb{F}_q$, then the group of $m$-torsion points
on an elliptic curve $E / \mathbb{F}_q$, denoted $E[m]$, is isomorphic to 
$(\mathbb{Z}/m\mathbb{Z}) \times (\mathbb{Z}/m \mathbb{Z})$. Fix an elliptic curve $E /\mathbb{F}_q$
and a prime $\ell$ that is not the characteristic of $\mathbb{F}_q$.
Let $P$ and $Q$ generate the group $E[\ell]$. 
A {\em distortion map} on $E$ is an endomorphism $\phi$ of $E$
such that $\phi(P) \notin \langle P \rangle$. A distortion map can be used to solve the DDH problem on the group $\langle P \rangle$ as follows: Given
a triple $R,S,T$ of points belonging to the group generated by $P$, we check whether $\mathbf{e}_{\ell}(R,\phi(S)) = \mathbf{e}_{\ell}(P,\phi(T))$, 
where $\mathbf{e}_{\ell}$ is the Weil pairing on the $\ell$-torsion points. It follows from well known
properties of the Weil pairing that this check succeeds if and only if $R = aP$, $S = bP$ and $T = abP$.
Under the assumptions that $P$ and $Q$ are both defined over $\mathbb{F}_{q^k}$, where $k$ is not large (say,
bounded by a fixed polynomial in $\log(q)$), and that $\phi$ can be computed in polynomial time, the DDH problem
can be solved in polynomial time using this idea. 
If $P$ and $Q$ are not eigenvectors for the Frobenius map, then in many cases
one can use the trace map as a distortion map (see \cite{gr04}). 
For this reason, we will concentrate only on the subgroups that are Frobenius eigenspaces.\\

It is known that distortion maps exist on supersingular elliptic curves (\cite{ver01,gr04}), 
and that distortion maps that do not commute with the Frobenius do not exist on ordinary elliptic curves (see \cite{ver01}
or \cite{ver04} Theorem 6). The latter implies that distortion maps
{\em do not} exist for ordinary elliptic curves with embedding degree $> 1$. The embedding degree, (say) $k$, is the order of $q$ 
in the group $\left( \mathbb{Z}/\ell\mathbb{Z}\right)^*$. A theorem of Balasubramanian and Koblitz (\cite{bk98} Theorem 1) says that if $E(\mathbb{F}_q)$ contains an $\ell$-torsion point and $k > 1$, then $E[\ell] \subseteq \mathbb{F}_{q^k}$.
Thus, the only remaining cases where the existence of Distortion maps is not known are the cases when the
embedding degree $k$ is $1$. If the embedding degree is $1$ and $E(\mathbb{F}_q)$ contains an $\ell$-torsion point, then
there are two possibilities: either $E[\ell](\mathbb{F}_q)$ is cyclic or $E[\ell] \subseteq E(\mathbb{F}_q)$. In the former situation there are no distortion maps (by \cite{ver04} Theorem 6). However, the Tate pairing
can be used to solve DDH efficiently in this case (see the comments in \cite{gr04} following Remark 2.2). Thus, the only 
case in which the question of the existence of a distortion map remains open is when $E[\ell] \subseteq E(\mathbb{F}_q)$. 
In this article we characterize the existence of distortion maps for this case.

\section{The Proof}
Let $k$ be a finite field, $\mathbb{F}_q \supseteq k$ and $E/k$ be an ordinary elliptic curve. Suppose $\ell$ is a prime
such that $E[\ell] \subseteq \mathbb{F}_q$ but no point of exact order $\ell$
is defined over a smaller field.\medskip

To study the existence of distortion maps, we study the reduction of the ring $\End(E)$ modulo $\ell$. 
Our principal tool is the following observation: 
If $\alpha \in \End(E)$ has field polynomial $f(x) \in \mathbb{Z}[x]$, then $f ~\rm{mod}~ \ell$
is the characteristic equation of the action of $\alpha$ on $E[\ell]$. \medskip

Let $\pi$ be the $q$-th power Frobenius endomorphism on $E$ and let $\phi^2 - t \phi + q = 0$ be its
characteristic equation. We know that $t \equiv 2 \mod \ell$ and $q \equiv 1 \mod \ell$
as the full $\ell$-torsion is defined over $\mathbb{F}_q$. \\

Let $\mathcal{O} = \mathrm{End}(E)$, $K = \mathcal{O} \otimes \mathbb{Q}$ and $\mathcal{O}_K$
the maximal order in $K$. We have the inclusions $\mathbb{Z}[\pi] \subseteq \mathcal{O} \subseteq \mathcal{O}_K$.
Since $t^2 - 4q = 0 \mod \ell$ we have that $\ell$ divides the product $[\mathcal{O} : \mathbb{Z}[\pi]][\mathcal{O}_K : \mathcal{O}] \mathrm{Disc}(K)$. 
The existence of distortion maps splits into cases depending on whether $\ell | [\mathcal{O}_K : \mathcal{O}]$ or $\ell | \mathrm{Disc}(K)$. Indeed, if $\ell | [\mathcal{O}_K : \mathcal{O}]$ there are no distortion maps, since the reduction modulo $\ell$ of every endomorphism is just multiplication by scalar. \\

In the following
we assume that $\ell \not|~ [\mathcal{O}_K : \mathcal{O}]$ so that the conductor of $\mathcal{O}$ is prime to $\ell$. Under this
assumption we have that the residue class rings
\begin{align*}
\mathcal{O}_K/(\ell) \cong \mathcal{O}/(\ell).
\end{align*}

Suppose that $\ell {\not|}~ \mathrm{Disc}(K)$ and that $\ell$ is {\em inert} in $\mathcal{O}_K$, then $\mathcal{O}/(\ell) \cong \mathbb{F}_{\ell^2}$. Let $\alpha \in \mathcal{O}$ be an endomorphism such that $\alpha \mod (\ell)$ does not lie in $\mathbb{F}_{\ell}$. Then
the action of $\alpha$ on $E[\ell]$ is irreducible since its characteristic equation is irreducible over $\mathbb{F}_{\ell}$. 
Now $\alpha$ gives us a distortion map on $E[\ell]$ since no subgroup of order $\ell$ of $E[\ell]$ is stabilized by $\alpha$.\\

Now if $\ell {\not|}~ \mathrm{Disc}(K)$ and $\ell$ is {\em split} in $\mathcal{O}_K$, then $\mathcal{O}/(\ell) \cong \mathbb{F}_{\ell}[X]/(X-a)(X-b) \cong (\mathbb{Z}/\ell\mathbb{Z})^2$ (where $a \neq b$). The action of any $\alpha \in O_K$, that
corresponds to the image of $X$ in $\mathbb{F}_{\ell}[X]/(X-a)(X-b)$ under the isomorphism, is conjugate to $\begin{pmatrix} \gamma & 0 \\
0 & \delta\end{pmatrix}$. Thus, distortion maps exist for all but two of the subgroups of $E[\ell]$.\\

Suppose that $\ell | \mathrm{Disc}(K)$ so that $\ell$ is {\em ramified} in $\mathcal{O}_K$, then $\mathcal{O}/(\ell) \cong  \mathbb{F_{\ell}}[X]/(X-a)^2$. Consider
the map $\alpha \in \mathcal{O}$ that corresponds to the image of $X$ in the ring $\mathbb{F}_{\ell}[X]/(X-a)^2$.
The action of $\alpha$
on $E[\ell]$ is conjugate to $\begin{pmatrix} 1 & \beta \\ 0 & 1 \end{pmatrix}$.
Note that $\beta \neq 0$, for if $\beta = 0$ then $\mathcal{O}/(\ell) \cong \mathbb{Z}/\ell \mathbb{Z}$,
but we know that $\mathcal{O}$ is rank $2$ over $\mathbb{Z}/\ell\mathbb{Z}$ since $\ell$ is ramified in $\mathcal{O}_K$ and does not divide the
conductor of $\mathcal{O}$. Thus, distortion maps exist for all but one subgroup of $E[\ell]$.\\

In summary, we have:
\begin{Thm}\label{thm_main:this}
Let $k$ be a finite field, $\mathbb{F}_q \supseteq k$ and $E/k$ be an ordinary elliptic curve
whose endomorphism ring is $\mathcal{O}$, an order in an imaginary quadratic field $\mathcal{O}$. Suppose $\ell$ is a prime
such that $E[\ell] \subseteq \mathbb{F}_q$ but no point of exact order $\ell$
is defined over a smaller field.\\
\begin{enumerate}
\item If $\ell ~|~ [\mathcal{O}_K : \mathcal{O}]$ there are no distortion maps.
\item If $\ell {\not|}~[\mathcal{O}_K : \mathcal{O}]\mathrm{Disc}(K)$ and
	\begin{enumerate}
	\item $\ell$ is inert in $\mathcal{O}_K$, then there are distortion maps for every (order $\ell$) subgroup of $E[\ell]$;
	\item $\ell$ is split in $\mathcal{O}_K$, then all but two subgroups of $E[\ell]$ have distortion maps.
	\end{enumerate}
\item If $\ell {\not|}~[\mathcal{O}_K:\mathcal{O}]$
and $\ell ~|~ \mathrm{Disc}(K)$ so that $\ell$ is ramified in $\mathcal{O}_K$, then all (except one) subgroups
of $E[\ell]$ have distortion maps.
\end{enumerate}
\end{Thm}

\section{Examples}
In this section, we give examples to illustrate that all the cases in Theorem \ref{thm_main:this} do occur.

\begin{Exa} Consider the elliptic curve $E: y^2 = x^3 + x$ over $\mathbb{Q}$. $E$ has complex multiplication
by $\mathbb{Z}[\imath]$ and has good reduction at all odd primes. Let $p$ be a prime such that $p \equiv 1 \mod 4$,
$\tilde{E}$ be the reduction of $E$ modulo $p$, and let $\imath^2 = -1 \mod p$.
Then $\tilde{E}[2] \subseteq \tilde{E}(\mathbb{F}_p)$ and $\tilde{E}[2]$ is 
$\{0_{\tilde{E}},(0,0),(\imath,0),(-\imath,0)\}$ where $0_{\tilde{E}}$
is the identity element. The map $[\imath]$ is an endomorphism that sends $(x,y) \mapsto (- x, \imath y)$.
It is easy to see that the map $[\imath]$ preserves the subgroup $\langle (0,0) \rangle$
and interchanges the remaining two subgroups, of order $2$, of $\tilde{E}[2]$. Note, that Deuring's reduction
theorem tells us that $\End(\tilde{E}) \cong \mathbb{Z}[i]$. Furthermore, in this case the subring $\mathbb{Z}[\pi]$
generated by the Frobenius is usually a smaller ring. Indeed, if $t$ is the trace of Frobenius
and $t^2 - 4p = -4b^2$, then the conductor of the order $\mathbb{Z}[\pi]$ is $b$. Now $b$ is at least $2$, since $t \equiv 2 \mod 4$, 
so $(t/2)$ is odd and we must have $p = (t/2)^2 + b^2$. 
Thus, case (3) of Theorem \ref{thm_main:this} applies and matches with what we observe for the $2$-torsion.
\end{Exa}

\begin{Exa}\label{exa_inert:this}(Suggested by anonymous reviewer). Let $E$ be the curve over $\mathbb{F}_{701}$ 
given by the equation $y^2 = x^3 - 35x + 98$.
Then $\End(E) = \mathbb{Z}[\frac{1+\sqrt{-7}}{2}]$ which is the maximal order in $\mathbb{Q}(\sqrt{-7})$. 
The order $\mathbb{Z}[\pi]$ has conductor $10$ in $\End(E)$.
The $5$-torsion is $\mathbb{F}_{701}$ rational, and moreover, $5$ is inert in $\End(E)$. Theorem \ref{thm_main:this} (2a)
shows that every subgroup of $E[5]$ admits a distortion map. Indeed, the map corresponding
to multiplication by $\alpha = \frac{1 + \sqrt{-7}}{2}$ is given by (\cite{sil94} Chapter II, Proposition 2.3.1 (iii))
\begin{align*}
[\alpha](x,y) = \left( 
	\alpha^{-2}\left( x - \frac{7(1-\alpha)^4}{x+\alpha^2 - 2} \right),
	\alpha^{-3}y\left( 1 + \frac{7(1-\alpha)^4}{(x+\alpha^2-2)^2}\right)
	\right).
\end{align*}

Let us check this for the group generated by the $5$-torsion point $P$ (with affine coordinates) $P = (224, 31)$.
Since $\alpha = 386 \in \mathbb{F}_{701}$, this tells us that $[\alpha](P) = (173,194)$. One checks
that the Weil pairing $\mathbf{e}_{5}(P,[\alpha](P)) = 464 \neq 1$. Thus, $[\alpha]$ works as a distortion map for
the group generated by $P$. \\

Now the $5$-torsion of $E$ is generated by $P$ and the point $Q = (573,450)$. A similar
computation shows that $[\alpha](Q) = (463,495)$. Also, $\mathbf{e}_5(Q,[\alpha]Q) = 89 \neq 1$. Again, this shows
that $[\alpha]$ works as a distortion map. \\

Given these calculations it is not hard to find the matrix of the action of $[\alpha]$ on $E[5]$ relative to
the basis $P,Q$ 
\begin{align*}
[\alpha] = \begin{pmatrix} 0 & -1 \\ 2 & 1 \end{pmatrix}.
\end{align*}
The characteristic polynomial of this matrix is irreducible modulo $5$ and thus the action on $E[5]$ is irreducible.
\end{Exa}

\begin{Exa} One can use the elliptic curve $E$ from Example \ref{exa_inert:this} to illustrate case (2b) of Theorem \ref{thm_main:this}.
This time we look at $E[2]$ (also contained in $\mathbb{F}_{701}$)
which is generated by the points $P = (319,0)$ and $Q = (389,0)$. The prime $2$ splits completely
in $\End(E)$. The proof of Theorem \ref{thm_main:this} tells us that the characteristic
polynomial of the action of the endomorphism $[\alpha]$ has two
distinct roots and would work as a distortion map for all but two subgroups of $E[2]$. Now the minimal polynomial
$\alpha$ is $x^2 - x + 2$ and modulo $2$ this splits as $x(x+1)$. Thus the action of $[\alpha]$ on $E[2]$ will have 
two eigenvectors, with eigenvalues $0$ and $1$ respectively. It is easy to check given the formula for $[\alpha]$
that indeed $[\alpha](P) = 0_E$ and $[\alpha](Q) = Q$.
\end{Exa}

\begin{Exa}
In this example we illustrate that case (1) of Theorem \ref{thm_main:this} also occurs.
Consider the curve $E / \mathbb{Q}$ given by the Weierstrass equation
\begin{align*}
y^2 = x^3 - \frac{3375}{121}x + \frac{6750}{121}.
\end{align*}
The $j$-invariant of E is $2^4 3^3 5^3$
and the conductor of $E$ is $108900$. $E$ has CM by the order of conductor $2$ in $\mathbb{Q}(\sqrt{-3})$.
Thus $\End(E) \cong \mathbb{Z} + 2 \mathcal{O}_K$ where $\mathcal{O}_K = \mathbb{Z} + \frac{1}{2}(1 + \sqrt{-3})\mathbb{Z}$.
$E$ has good reduction at the prime $13$ and one sees that the reduction $\tilde{E}$ has $\mathbb{F}_{13}$-rational
$2$-torsion. Now $\End(\tilde{E}) \cong \End(E)$ by the Deuring reduction theorem (\cite{lan87} Chapter 13 \S4, Theorem 12),
but $\End(\tilde{E}) \mod 2 \cong (\mathbb{Z}/2\mathbb{Z})$ and so there are no distortion maps.
\end{Exa}

\end{document}